\documentclass[preprint]{ifacconf}

\usepackage{graphicx}      % include this line if your document contains figures
\usepackage{natbib}        % required for bibliography
\usepackage{amsfonts}
\usepackage{amssymb}  % assumes amsmath package installed
\usepackage{amsmath}
\usepackage{bm}
\usepackage{color}
\usepackage{tikz}
\usepackage{pgfplots}
\usepackage{float}
\usepackage{mathtools}

\pgfplotsset{compat=newest}
\usetikzlibrary{plotmarks}
\usetikzlibrary{arrows.meta}
\usepgfplotslibrary{patchplots}
\newlength\fwidth
\newlength\hwidth

\newtheorem{remark}[thm]{Remark}
\newtheorem{assumption}[thm]{Assumption}

\newfloat{algorithm}{t}{}

\newcommand{\norm}[1]{\left\lVert#1\right\rVert}
\newcommand\numberthis{\addtocounter{equation}{1}\tag{\theequation}}
\renewcommand{\refeq}[1]{\stackrel{\mathclap{#1}}{=}}
\newcommand{\refleq}[1]{\stackrel{\mathclap{#1}}{\leq}}

\graphicspath{{figures/}}

\begin{document}
\begin{frontmatter}
	
\vspace{-10ex}
	
{\textcopyright~2020 the authors. This work has been accepted to IFAC for publication under a Creative Commons Licence CC-BY-NC-ND}

\title{Online Gradient Descent for Linear Dynamical Systems} %\thanksref{footnoteinfo}} 

%\thanks[footnoteinfo]{Sponsor and financial support acknowledgment
%goes here. Paper titles should be written in uppercase and lowercase
%letters, not all uppercase.}

\author[First]{Marko Nonhoff} 
\author[First]{Matthias A. M\"uller}

\address[First]{Institute of Automatic Control, Leibniz University Hannover, 30167 Hannover, Germany.}

\begin{abstract}                % Abstract of not more than 250 words.
In this paper, online convex optimization is applied to the problem of controlling linear dynamical systems. An algorithm similar to online gradient descent, which can handle time-varying and unknown cost functions, is proposed. Then, performance guarantees are derived in terms of regret analysis. We show that the proposed control scheme achieves sublinear regret if the variation of the cost functions is sublinear. In addition, as a special case, the system converges to the optimal equilibrium if the cost functions are invariant after some finite time. Finally, the performance of the resulting closed loop is illustrated by numerical simulations.
\end{abstract}


\end{frontmatter}
%===============================================================================

\section{Introduction}

Online convex optimization is an extension of classical numerical optimization to the case where an algorithm operates online in an unknown environment. Whereas in convex optimization the goal is to minimize a given cost function subject to known constraints \citep{Boyd2004,Nesterov2018}, in an online convex optimization (OCO) problem the cost function to be minimized is time-varying and the algorithm only has access to past information. Specifically, at every time $t$, the algorithm has to choose an action $y_t \in \mathbb{Y}$ from an action set $\mathbb{Y}$ based on the actions chosen at previous time instances and the corresponding observed cost functions. Then, the environment reveals a new cost function $L_t: \mathbb{Y} \rightarrow \mathbb{R}$, which leads to the cost $L_t(y_t)$. The goal is to minimize the total cost in $T$ stages. The classical OCO framework was introduced in \citep{Zinkevich2003} and has received considerable interest as a tool for online optimization and learning (see \cite{Shwartz2012,Hazan2016} for an overview). The performance of OCO algorithms is commonly characterized by regret, which is defined as the gap between the algorithm's performance and some offline optimum in hindsight. In \citep{Hazan2007}, several algorithms which achieve low static regret, i.e., low regret with respect to the best constant action, are presented. Dynamic benchmarks are proposed in \citep{Besbes2015, Mokhtari2016}. Sublinear regret is generally desirable because it implies that the algorithm's performance is asymptotically on average no worse than the benchmark. One major advantage of the OCO framework is its ability to handle time-invariant as well as time-varying constraints \citep{Paternain2016,Cao2019}.

In particular, an online version of gradient descent termed online gradient descent (OGD) has proven to be a simple algorithm that achieves low regret in the OCO setting \citep{Hazan2016}. In OGD, at every time instant $t$, the action $y_t$ is chosen as $y_t = \Pi_\mathbb{Y}(y_{t-1} - \gamma \nabla f_{t-1}(y_{t-1}))$, where $\Pi_\mathbb{Y}(y)$ denotes a projection of a point $y$ onto the convex constraint set $\mathbb{Y}$, $\gamma \in \mathbb{R}$ is a step size parameter, and $f_t(y)$ is the cost function to be minimized. Hence, instead of solving an optimization problem at every time instant $t$, only one gradient descent step on the previous cost function is employed to reduce computational complexity.

Whereas classical OCO does not consider coupling \linebreak between time instances and, therefore, no underlying \linebreak dynamical system, some combinations of OCO with dynamical models have been studied recently. On the one hand, in \citep{Hall2015}, algorithms for online prediction of incoming data are proposed, where the data is generated by a dynamical system. It is shown that the regret of the proposed algorithms is low if the environment follows the model of the underlying system. On the other hand, coupling between time instances in the OCO setting has been considered by introducing a switching or cost $d(y_t-y_{t-1})$ to study the effect of a time coupled cost function \citep{Tanaka2006,Lin2013}. In \citep{Li2018}, the switching cost is chosen as $d(y_t-y_{t-1})=\frac{\beta}{2}\norm{y_t - y_{t-1}}^2$, where $\beta \in \mathbb{R}$ is a weighting parameter. This can be interpreted as an additional quadratic cost on the input $u_{t-1}$ of a single integrator system $y_t = y_{t-1} + u_{t-1}$. This approach is extended to general linear systems in \citep{Li2019}. However, this work strongly focuses on the case where predictions of future cost functions are available to the algorithm. It is shown that the algorithm's regret can be reduced substantially by utilizing predictions. Moreover, in \citep{Abbasi2014,Cohen2018,Akbari2019}, linear dynamical systems and quadratic cost functions are considered. Therein, the best linear controller is chosen as the benchmark in the definition of regret and the proposed algorithms apply OCO to optimize over the set of stable, linear policies. This method is generalized to general convex cost functions in \citep{Agarwal2019}. A different approach is taken in \citep{Colombino2020}, where online optimization is used to steer a linear dynamical system to the solutions of time-varying convex optimization problems.

In contrast, many control algorithms able to handle dynamical systems equipped with a cost function as well as state and input constraints exist. In particular, Model Predictive Control (MPC) is able to minimize a given cost function while taking constraints explicitly into account \citep{Rawlings2009}. However, classical MPC techniques require solving a potentially large optimization problem at every time instant. Therefore, inexact or suboptimal MPC algorithms have been proposed, which only complete a finite number of optimization iterations at every time instant \citep{Scokaert1999, Diehl2005}. %Optimization algorithms in inexact MPC schemes include, for example, interior point \citep{Pavlov2019}, active set \citep{Zeilinger2011}, and (fast) gradient methods \citep{Richter2011, Patrinos2014}, for which stability, feasibility, and computational complexity have been studied. However, there are limited results on the closed-loop performance of suboptimal MPC.
Research in the field of suboptimal MPC typically focuses on stability, feasibility, and computational complexity. However, there are limited results on the closed-loop performance of suboptimal MPC schemes.

In this work, OCO for controlling a general linear controllable dynamical system equipped with a general convex cost function $L_t(x,u)$ is considered. We focus on the cases where the cost functions are strongly convex and smooth, and we propose an algorithm which is similar to online gradient descent in classical OCO. We show that the proposed algorithm attains sublinear regret if the variation of the cost functions is sublinear in time. In contrast to existing results in OCO, which do not consider an underlying dynamical system, this requires novel algorithm design and analysis techniques since the dynamical system cannot move in an arbitrary direction in a single time instant. Moreover, we do not restrict the proposed algorithm to the set of stable, linear controllers. Compared to suboptimal MPC, we explicitly consider a time-varying stage cost which is usually not possible in the MPC literature. In particular, in MPC the cost function needs to be known a priori in order to obtain performance guarantees.

This paper is organized as follows. Section~\ref{sec:setting} defines the problem setting and states our proposed algorithm. A theoretical analysis of the closed loop and the main theorem are given in Section~\ref{sec:thmproof}. In Section~\ref{sec:sim}, we illustrate the algorithm's performance by numerical simulations. Section~\ref{sec:conclusion} concludes the paper.

\textit{Notation}: For a vector $x \in \mathbb{R}^n$, $\norm{x}$ denotes the Euclidean norm, whereas for a matrix $A \in \mathbb{R}^{n\times m}$, $\norm{A}$ denotes the corresponding induced matrix norm and $A^T$ the transposed of the matrix $A$. We define by $\mathbb{N}_{[a,b]}$ the set of natural numbers in the interval $[a,b]$. The gradient of a function $f(x)$ is denoted by $\nabla f(x)$. Additionally, $\bm{I}_n$ is the identity matrix of size $n \times n$, and $\bm 0_{n \times m}$ is the matrix of all zeros of size $n \times m$.

\section{Setting and Algorithm} \label{sec:setting}
\subsection{Problem setup}

We consider discrete-time linear systems of the form
\begin{align}
	x_{t} &= Ax_{t-1} + Bu_t \label{eq:sysdyn}
\end{align}
with a given initial condition $x_0 \in \mathbb{R}^n$, where $x_t \in \mathbb{R}^n$ are the states of the system, $u_t \in \mathbb{R}^m$ are the control inputs, and $A \in \mathbb{R}^{n \times n}$, $B \in \mathbb{R}^{n \times m}$.
\begin{remark}
	Note that we slightly deviate from the usual notation for linear systems (i.e., $x_{t+1} = Ax_t + Bu_t$) to facilitate notation in the proposed OCO approach.
\end{remark}
At every time instant $t \in \mathbb{N}_{[1,T]}$, we choose a control action $u_t$ which is applied to system \eqref{eq:sysdyn}. Then, afterwards, a cost function $L_{t}: \mathbb{R}^n \times \mathbb{R}^m \rightarrow \mathbb{R}$ is revealed which results in the cost $L_{t}(x_{t},u_t)$ before we move on to the next time step. As is standard in OCO, we measure our algorithm's performance by regret. In our case, regret is defined as
\begin{align*}
	\mathcal{R} &:= \sum_{t=1}^T L_{t}(x_{t},u_t) - L_{t}(x^*_{t},u^*_t).
\end{align*}
Here, the state and input sequences $\bm x^* = (x_1^*, \dots, x_T^*)$ and $\bm u^* = (u_1^*, \dots, u_T^*)$ are defined as the solution to the optimization problem
\begin{align*}
	\min_{u_1,x_1,\dots,u_T,x_T} &\sum_{t=1}^T L_t(x_t,u_t)
	&&\text{s.t.} &&&x_t = Ax_{t-1} + Bu_t.
\end{align*}
Hence, $(x^*_t,u^*_t)$ denote the optimal states and inputs at time $t$ in hindsight, with full knowledge of the cost functions $L_t$. The regret $\mathcal{R}$ can therefore be viewed as a measure of how much we regret receiving information about the cost functions $L_t$ only after we choose a control input $u_{t}$. We do not consider any cost on the initial condition $x_0$ at time $t=0$ since it cannot be influenced by the algorithm's decisions, i.e., control inputs $u_t$.

Similar to various works in OCO (see, e.g., \cite{Mokhtari2016,Li2018}), we assume the cost functions $L_t$ to be separable, strongly convex, and smooth as stated in the following assumption.
\begin{assumption} \label{assump:CostFunctionProp}
	For every $t \in \mathbb{N}_{[0,T]}$, the cost function $L_t$ satisfies the following conditions:
	\begin{enumerate}
		\item $L_t(x,u) = f_t^x(x) + f_t^u(u)$,
		\item $f_t^x(x)$ is $\alpha_x$-strongly convex and $l_x$-smooth,
		\item $f_t^u(u)$ is $\alpha_x$-strongly convex and $l_x$-smooth,
	\end{enumerate}	
\end{assumption}
where an $\alpha$-strongly convex function $f: \mathbb{R}^n \rightarrow \mathbb{R}$ satisfies 
\begin{equation*}
f(y) \geq f(x) + \left\langle \nabla f(x),y-x \right\rangle + \frac{\alpha}{2} \norm{y-x}^2,~\forall x,y \in \mathbb{R}^n,
\end{equation*}
and $l$-smoothness means that $f$ satisfies
\begin{equation*}
f(y) \leq f(x) + \left\langle \nabla f(x),y-x \right\rangle + \frac{l}{2} \norm{y-x}^2,~\forall x,y \in \mathbb{R}^n.
\end{equation*}
We note the fact that $l$-smoothness of a convex function $f$ implies the following Lipschitz condition on the gradient of $f$ \citep{Nesterov2018} $\norm{\nabla f(x) - \nabla f(y)} \leq l \norm{x-y}$.

Moreover, we define \mbox{$\theta_t = \arg \min_x f_t^x(x)$} and \linebreak \mbox{$\eta_t = \arg \min_u f_t^u(u)$}. In the following, we assume that the minima are attained and therefore finite and, due to convexity, unique. Hence, at each time $t$ when the cost $L_t(x,u)$ is measured, that is $t \in \mathbb{N}_{[1,T]}$, $(\theta_t, \eta_t)$ is the minimizer of $L_t(x,u)$. In contrast to the trajectories $\bm x^*$ and $\bm u^*$, the sequences $\bm \theta = (\theta_1, \dots, \theta_T)$ and $\bm \eta = (\eta_1,\dots,\eta_T)$ do in general \emph{not} satisfy the system dynamics \eqref{eq:sysdyn}.

So far, the optimizer $(\theta_t, \eta_t)$ and the cost functions $L_t$ are only defined for $t \in \mathbb{N}_{[1,T]}$. For the remainder of this work, we fix without loss of generality $L_0(x,u)=f_0^x(x)+f_0^u(u)$ such that Assumption \ref{assump:CostFunctionProp} is satisfied, $\arg \min_x f_0^x(x) = \theta_0$, and $\arg \min_u f_0^u(u) = \eta_0$. The values of $\theta_0$ and $\eta_0$ can be defined arbitrarily, and a convenient choice for our subsequent analysis is given below.

If the cost function $L_t$ is allowed to change arbitrarily at every time step, we cannot expect to achieve a low regret. Therefore, we consider the path length as a measure for the variation of the cost functions $L_t$. Path length is defined as the accumulative absolute difference of the optimizer of the cost function at two consecutive time steps:
\begin{align*}
\text{Path length} := \sum_{t=1}^T \norm{\theta_t - \theta_{t-1}} + \sum_{t=1}^T \norm{\eta_t - \eta_{t-1}}.
\end{align*}
Path length plays an important role in dynamic regret analysis \citep{Mokhtari2016}. In \citep{Li2018} it is shown for single integrator systems that sublinear regret can be achieved if the path length is sublinear in $T$.

In this paper we consider tracking cost functions as \linebreak described in the following assumption.
\begin{assumption} \label{assump:TrackingCostFcn} For all $t \in \mathbb{N}_{[1,T]}$, $\theta_t$ and $\eta_t$ satisfy
	\begin{align*}
		\theta_t = A\theta_t + B\eta_t.
	\end{align*}
\end{assumption}
Assumption \ref{assump:TrackingCostFcn} states that the minimum $(\theta_t,\eta_t)$ of the cost function $L_t(x,u)$ is a steady-state with respect to the system dynamics \eqref{eq:sysdyn}, meaning that the control objective is to track a priori unknown and online changing setpoints. Relaxing this assumption to general convex cost functions (termed \emph{economic} cost functions in the context of MPC \citep{Faulwasser2018}) is part of our ongoing work.

Last, we assume that system \eqref{eq:sysdyn} is controllable and require the norm of $A \in \mathbb{R}^{n \times n}$ to be bounded as follows:
\begin{assumption} \label{assump:SysCond}
	$(A,B)$ is controllable with controllability index $\mu$, i.e.,
	\begin{align*}
	\text{rank}~S_c = \text{rank}~(B~AB~\ldots~A^{\mu-1}B) &= n.
	\end{align*}
	Moreover, $A$ satisfies $\norm{A} < \frac{l_x + \alpha_x}{2(l_x - \alpha_x)}$.
\end{assumption}

The norm of $A$ can be interpreted as a measure for the stability/instability of system~\eqref{eq:sysdyn}. Achieving low regret for a sequence of cost functions $f_t^x(x)$ by applying an algorithm similar to OGD means that the gradient descent step needs to counteract the instability of the system. Hence, it is natural to require an upper bound for the norm of $A$ in terms of the smoothness and convexity of the stage cost $f_t^x(x)$. In particular, if the cost function is given by $L_t(x,u) = \frac{\beta}{2} \norm{x - \theta_t}^2 + f_t^u(u)$ for some $\beta > 0$, we obtain $\alpha_x = l_x = \beta$. Hence, in this case, any controllable system satisfies Assumption \ref{assump:SysCond}.

\subsection{Online gradient descent for linear systems}

\relpenalty=9999
\binoppenalty=9999

Before we state our algorithm, we first define two useful matrices. The matrix $W = \begin{pmatrix} \bm{0}_{m \times (\mu-1)m} & \bm{0}_{m \times m} \\ I_{(\mu-1)m} & \bm{0}_{(\mu-1)m \times m} \end{pmatrix}$ shifts a vector by $m$ components, whereas the matrix $e = \begin{pmatrix} \bm{0}_{m \times (\mu-1)m} & I_m \end{pmatrix}$ extracts the last $m$ components.

\begin{algorithm} \centering
	{
		\setlength\belowdisplayskip{0pt}
		\setlength\abovedisplayskip{0pt}
		\fbox{\parbox{.9\linewidth}{
				
				\textbf{Algorithm 1}~(OGD for linear dynamical systems)
				
				\vspace{-6pt}						
				\rule{.44\textwidth}{.5pt}
				
				Given step sizes $\gamma_v$ and $\gamma_x$, initialization $v_0$, $x_0$, and state vector $x_{t-1}$. At time $t \in [1,T]$:
				\begin{align}
				\intertext{Input OGD}
				v_t &= v_{t-1} - \gamma_v \nabla f_{t-1}^u(v_{t-1}) \label{eq:GenB_InputOGD} \\
				\intertext{Prediction} \hat V_t &= \begin{pmatrix} v_t \\ \vdots \\ v_t \end{pmatrix} \in \mathbb{R}^{\mu m} \label{eq:GenB_InputPrediction} \\
				\hat x_{t+\mu-1} &= A^{\mu} x_{t-1} + S_c \hat V_t + S_c \sum_{i=1}^{\mu-1} W^i g_{t-i} \label{eq:GenB_StatePrediction} \\
				\intertext{State OGD} g_t &= -\gamma_x S_c^T \left(S_c S_c^T \right)^{-1} \nabla f_{t-1}^x(\hat x_{t+\mu-1}) \label{eq:GenB_StateOGD} \\
				\intertext{Output}
				u_t &= v_t + \sum_{i=0}^{\mu - 1} e W^i g_{t-i} \label{eq:GenB_Output}
				\end{align}
				\vspace{-6pt}
		}}
	}
\end{algorithm}

The proposed OCO scheme is given in Algorithm~1, where we set $g_t = 0$ if $t < 1$ in \eqref{eq:GenB_StatePrediction} and \eqref{eq:GenB_Output}. Note that the inverse $(S_cS_c^T)^{-1}$ in \eqref{eq:GenB_StateOGD} exists due to controllability in Assumption~\ref{assump:SysCond}. In our setting, at every time step $t$, given the previous state $x_{t-1}$ and cost function $L_{t-1}(x,u)$, Algorithm~1 computes a control input $u_t$ which is then applied to system \eqref{eq:sysdyn}. Afterwards, a new cost function $L_t$ is revealed resulting in the cost $L_t(x_t,u_t)$. 

Since no cost function is known at the first time instant, a standard method in OCO is to apply an arbitrary initialization input $v_0$. At time $t=1$, Algorithm~1 computes $v_1 = v_0 - \gamma_v\nabla f_0^u(v_0)$ and $g_1 = -\gamma_x S_c^T(S_cS_c^T)^{-1}\nabla f_0^x(\hat x_\mu)$. Therefore, we fix $\theta_0 = \hat x_\mu$ and $\eta_0 = v_0$, which yields $v_1 = v_0$, $g_1 = 0$, and, hence, $u_1 = v_1 = v_0$.

\begin{figure}
	\centering \small
	\def\svgwidth{.4\textwidth}
	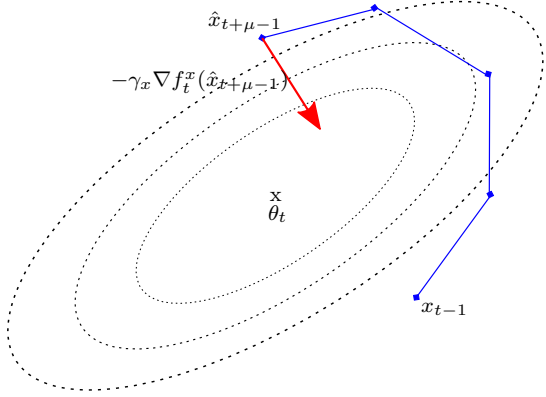
	\caption{Illustration of Algorithm~1. Dashed: contour lines of $L_t(x,u)$; blue: Predicted system states for the next $\mu$ time steps; red: desired gradient descent step.}
	\label{fig:SchematicPrediction}
	\setlength{\belowdisplayskip}{5pt}
\end{figure}

Roughly speaking, the proposed algorithm employs OGD twice to seek the optimal input $\eta_t$ and the optimal state $\theta_t$. First, we apply OGD in~\eqref{eq:GenB_InputOGD} to track the optimal input. Next, we would like to apply OGD again on the states of the system which would yield $Bu_t = Bv_t - \gamma_x \nabla f_{t-1}(Ax_{t-1}+Bv_t)$. Unfortunately, this is not possible if the system is not \mbox{$1$-step} controllable, i.e., $\text{rank}(B) < n$. Instead, as illustrated in Figure~\ref{fig:SchematicPrediction}, the algorithm predicts an input sequence for the next $\mu$ time steps and the corresponding system state $\mu$ time steps ahead in \eqref{eq:GenB_StatePrediction}. Then, an additional input sequence $\bm u^{OGD} = (u_0^{OGD}, \dots, u_{\mu-1}^{OGD})$ for the next $\mu$ time instances is determined such that application of both computed input sequences results in a gradient descent step on the previous cost function $\mu$ time steps in the future. Hence, we require $S_c g_t = -\gamma_x \nabla f_{t-1}^x(\hat x_{t+\mu-1})$, where $g_t = \begin{pmatrix} (u_{\mu-1}^{OGD})^T & \dots & (u_{0}^{OGD})^T \end{pmatrix}^T \in \mathbb{R}^{\mu m}$ is the vector created by stacking the components of the additional input sequence. Moreover, because we want $v_t$ to converge to the optimal input $\eta_t$, we need $g_t$ such that it does not contribute much to the input cost $f_t^u(u_t)$. Therefore, we choose $g_t$ to be the solution of
\begin{align*}
	&\min_{g_t \in \mathbb{R}^{\mu m}} \norm{g_t}^2 &&\text{s.t. } &&&S_c g_t = -\gamma_x \nabla f_{t-1}^x(\hat x_{t+\mu-1}).
\end{align*}
Solving this optimization problem analytically yields \eqref{eq:GenB_StateOGD}. By employing the matrices $W$ and $e$ as in \eqref{eq:GenB_Output}, we extract the required input $u_{i}^{OGD}$ from $g_t$ at time $t+i,~i\in\mathbb{N}_{[0,\mu-1]}$, to complete one gradient descent step at time $t+\mu-1$.

\section{Theoretical Analysis} \label{sec:thmproof}

In this section, we state our main result which gives a bound on the regret of Algorithm~1. To shorten notation, let $\Theta_\tau = \sum_{t=1}^\tau \norm{\theta_t - \theta_{t-1}}^2$ and $H_\tau = \sum_{t=1}^\tau \norm{\eta_t - \eta_{t-1}}^2$.

\begin{samepage}
\begin{thm} \label{thm}
	Let Assumptions \ref{assump:CostFunctionProp}, \ref{assump:TrackingCostFcn}, and \ref{assump:SysCond} be satisfied and let $\alpha_u > \frac{2-\sqrt 2}{2+\sqrt2} l_u$. Given step sizes $\frac{2-\sqrt 2}{2\alpha_u} < \gamma_v \leq \frac{2}{l_u+\alpha_u}$ and $\frac{2\norm{A}-1}{2\norm{A}\alpha_x} < \gamma_x \leq \frac{2}{l_x + \alpha_x}$, there exist constants $\Lambda_\theta > 0$ and $\Lambda_\eta > 0$, such that for each $T \in \mathbb{N}_{\geq 1}$, the regret of Algorithm~1 can be upper bounded by
	\begin{equation*}
		\mathcal{R} \leq C_\mu + \Lambda_\theta \Theta_T + \Lambda_\eta H_T,
	\end{equation*}
	where $C_\mu = l_x/2 \sum_{t=1}^{\mu-1} \norm{x_t - \theta_t}^2$.
\end{thm}
\end{samepage}

The proof of Theorem \ref{thm} is given in the appendix.

Theorem \ref{thm} states that the regret of Algorithm~1 can be upper bounded, where the bound depends linearly on $\Theta_T$ and $H_T$, up to a constant $C_\mu$. First, note that the step sizes $\gamma_v$ and $\gamma_x$ are well defined. Due to the lower bound on $\alpha_u$ in Theorem~\ref{thm} we have \mbox{$\frac{2-\sqrt 2}{2\alpha_u} = \frac{2(2-\sqrt 2)}{(2+\sqrt 2) \alpha_u + (2 - \sqrt 2) \alpha_u} < \frac{2(2-\sqrt 2)}{(2-\sqrt 2)l_u + (2 - \sqrt 2) \alpha_u} = \frac{2}{l_u + \alpha_u}$}. In addition, the upper bound on $\norm{A}$ in Assumption \ref{assump:SysCond} yields $\frac{2\norm{A}-1}{2\norm{A}\alpha_x} = \frac{1}{\alpha_x} - \frac{1}{2\norm{A}\alpha_x} < \frac{1}{\alpha_x} - \frac{l_x-\alpha_x}{(l_x + \alpha_x) \alpha_x} = \frac{2}{l_x + \alpha_x}$. Second, $C_\mu = l_x/2 \sum_{t=1}^{\mu-1} \norm{x_t - \theta_t}^2$ can be bounded in terms of $\Theta_\mu$ and $H_\mu$ as well. While we omit a detailed derivation due to space limitations, we note that, since $\theta_t$ and $\eta_t$ are finite, $C_{\mu}$ is a constant (depending on $x_0$, $v_0$, $\theta_1$, $\dots$, $\theta_{\mu-1}$, and $\eta_1$, $\dots$, $\eta_{\mu-1}$) which is \emph{independent} of $T$. This is sufficient in order to attain sublinear regret (compare Corollary~\ref{cor} below). The constants $\Lambda_\theta$ and $\Lambda_\eta$, which are independent of $T$ as well, can be explicitly calculated as is shown in the proof. 

\begin{remark} \label{rem:regret}
	In the first step of bounding the regret in the proof of Theorem~\ref{thm}, we exploit optimality of $(\theta_t,\eta_t)$ to lower bound $f_t^x(x^*_t)$ and $f_t^u(u^*_t)$ by $f_t^x(\theta_t)$ and $f_t^u(\eta_t)$, respectively. Hence, we could also use the point-wise in time optima $\bm \theta$ and $\bm \eta$ as a benchmark in the definition of regret instead of the best possible trajectories $\bm x^*$ and $\bm u^*$ and would still achieve the same bound on the regret. Characterizing and exploiting properties of the trajectories $\bm x^*$ and $\bm u^*$ to achieve a less conservative regret bound is an interesting topic for future research. Furthermore, if $L_t(\theta_t,\eta_t) = 0$ for all $t \in \mathbb{N}_{[1,T]}$, the bound on the regret in Theorem~\ref{thm} is a bound on the total cost over $T$ stages.
\end{remark}

Having established an upper bound on the regret $\mathcal{R}$ of Algorithm~1, we now examine whether the regret is sublinear in $T$. The corollary below gives a sufficient condition for sublinear regret.

\begin{cor} \label{cor}
	Let $\theta_t \in \mathbb{D}_\theta$ and $\eta_t \in \mathbb{D}_\eta$ for all \mbox{$t \in \mathbb{N}$}, where $\mathbb{D}_\theta$ and $\mathbb{D}_\eta$ are compact sets. If the path length $\sum_{t=1}^T \norm{\theta_t - \theta_{t-1}} + \sum_{t=0}^T \norm{\eta_t-\eta_{t-1}}$ is sublinear in $T$, then the regret of Algorithm~1 is sublinear in $T$.
\end{cor}
\emph{Proof of Corollary \ref{cor}.} Compactness of $\mathbb{D}_\theta$ and $\mathbb{D}_\eta$ implies $\norm{\theta_t - \theta_{t-1}} \leq d_\theta$ for some $d_\theta \in \mathbb{R}$ and $\norm{\eta_t - \eta_{t-1}} \leq d_\eta$ for some $d_\eta \in \mathbb{R}$ for all $t \in \mathbb{N}$. Hence, we have $\Theta_T = \sum_{t=1}^T \norm{\theta_t - \theta_{t-1}}^2 \leq d_\theta \sum_{t=1}^T \norm{\theta_t - \theta_{t-1}}$, \linebreak $H_T = \sum_{t=1}^T \norm{\eta_t - \eta_{t-1}}^2 \leq d_\eta \sum_{t=1}^T \norm{\eta_t - \eta_{t-1}}$.
Thus, $\Theta_T$ and $H_T$ are sublinear in $T$. Theorem \ref{thm} states that the regret of Algorithm~1 is at most linear in $\Theta_T$ and $H_T$. Therefore, the regret of the proposed algorithm is sublinear in $T$. \hfill \hfill \qed

\begin{remark}
	Consider, as a special case, that for some finite time $t'$ the minimizer of the cost functions satisfy $\theta_{t} = \theta_{t'}$ and $\eta_t = \eta_{t'}$ for all $t \geq t'$. Optimality of $(\theta_t, \eta_t)$ and the fact that $\theta_t$ and $\eta_t$ are both finite implies
	\begin{align*}
		\mathcal{R} &\leq \sum_{t=1}^T f_t^x(x_t) - f_t^x(\theta_t) + f_t^u(u_t) - f_t^u(\eta_t) \\
		&\leq  C_\mu + \Lambda_\theta \Theta_T + \Lambda_\eta H_T = C_\mu + \Lambda_\theta \Theta_{t'} + \Lambda_\eta H_{t'} < \Lambda
	\end{align*}
	for some $\Lambda \in \mathbb{R}$ independent of $T$, where the second inequality is by Theorem~\ref{thm} (compare Remark \ref{rem:regret}). We also have \mbox{$f_t^x(x_t) - f_t^x(\theta_t) \geq 0$} as well as $f_t^u(u_t) - f_t^u(\eta_t) \geq 0$. Now, if we let $T \rightarrow \infty$, we obtain $\lim_{T \rightarrow \infty} x_t = \theta_{t'}$ and $\lim_{T \rightarrow \infty} u_t = \eta_{t'}$, i.e., the closed loop converges to the optimal equilibrium.
\end{remark}

\section{Simulations} \label{sec:sim}

We illustrate the effectiveness of the proposed algorithm through numerical simulations. We randomly choose the matrix $A=\begin{pmatrix} 1.05 & 0.7 & 1.75 \\ 0.35 & 0.7 & 1.05 \\ 1.4 &0.105 & 1.855 \end{pmatrix}$. We set $B=\begin{pmatrix} 1 & 0 & 1 \end{pmatrix}^T$ which yields a controllable system with $\mu = 3$. We set $T=30$, the cost function to $L_t(x,u) = f_t^x(x) + f_t^u(u) = \frac{1}{2}\norm{x-\theta_t}^2 + \frac{1}{2}\norm{u-\eta_t}^2$, and choose $\gamma_v = 0.98$ and $\gamma_x = 0.995$. It is easy to see that this choice satisfies all conditions in Theorem \ref{thm} since $f_t^x$ and $f_t^u$ are $\alpha$-strongly convex and $l$-smooth, where $\alpha = l = 1$. Algorithm~1 is initialized with $x_0 = \begin{pmatrix} 0 & 0 & 0 \end{pmatrix}^T$ and $v_0 = 0$.

\begin{figure}
	\centering \small
	\setlength\fwidth{.4\textwidth}
	\setlength\hwidth{.9\fwidth}
	\input{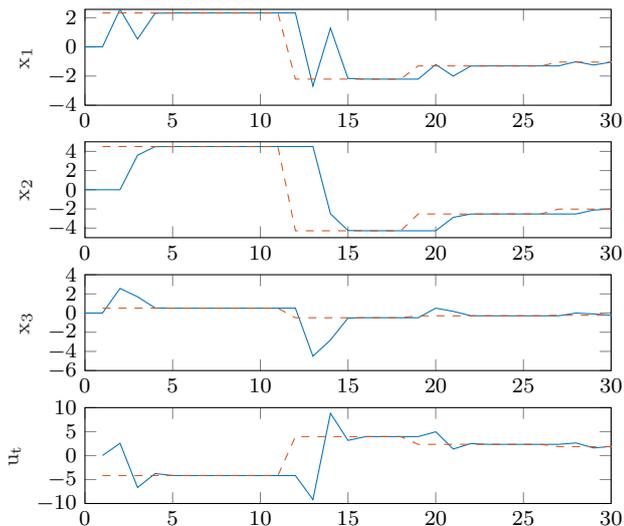}
	\vspace{-35pt}
	\caption{Blue, solid: The three states $x_i$, $i\in \{1,2,3\}$, and the input $u_t$; red, dashed: The optimal states $\theta_i$, $i \in \{1,2,3\}$, and input $\eta_t$.}
	\label{fig:StatesInputs}
\end{figure}

In the first experiment, the sequences $\bm \theta$ and $\bm \eta$ are chosen randomly, i.e., $\eta_1$ is sampled randomly from a uniform distribution over the interval $[-5, 5]$ and $\theta_1$ is calculated such that Assumption \ref{assump:TrackingCostFcn} is satisfied. Then, at every time instant $t \in \mathbb{N}_{[2,30]}$, the minimizer $(\theta_t,\eta_t)$ are modified with a probability of $0.1$. If they change, $\eta_t$ is again sampled randomly from the interval $[-5, 5]$ and $\theta_t$ is calculated accordingly. Hence, the control objective is to track a priori unknown setpoints. Figure~\ref{fig:StatesInputs} presents the resulting closed-loop trajectories for all three states and the input trajectory $u_t$. The closed loop closely tracks the desired setpoints, and it converges to the optimal states and control input within a few time steps whenever the cost function is changed.

\begin{figure} 
	\centering \small
		\setlength\fwidth{.4\textwidth}
		\setlength\hwidth{.7\fwidth}
		\input{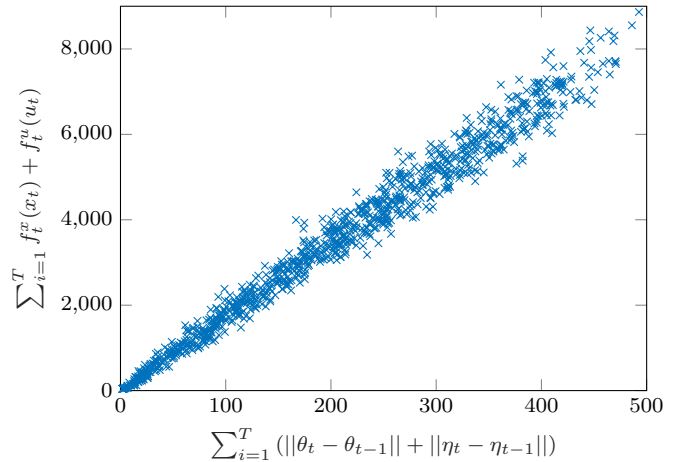}
		\vspace{-20pt}
		\caption{Total cost of $1000$ simulation runs over the path length of each run.}
		\label{fig:RandomSim}
\end{figure}

In a second experiment, the system and the cost function are chosen as before and $1000$ simulations are conducted, where we set $T=500$. At every time step $t$, the optimal state and input $(\theta_t,\eta_t)$ is changed with a probability of $0.25\frac{j}{1000}$, where $j \in \mathbb{N}_{[1,1000]}$ is the number of the simulation run. As before, if the cost function is changed, the new value for $\eta_t$ is sampled randomly from the interval $[-5,5]$ and $\theta_t$ is calculated such that Assumption \ref{assump:TrackingCostFcn} is satisfied. Thereby, we cover a wide range of path lengths. Figure \ref{fig:RandomSim} shows the resulting total cost and path length for each simulation run. Apparently, the total cost grows linearly with the path length as stated in Theorem~\ref{thm}.

\section{Conclusion} \label{sec:conclusion}

In this paper, we apply online convex optimization to linear dynamical systems equipped with a cost function and propose a first online algorithm for this problem. We derive a bound on the regret of the algorithm and show that it achieves sublinear regret if the variation of the cost functions, measured in terms of path length, is sublinear in time. The performance of the proposed algorithm is illustrated by numerical examples.

Since we do not consider any constraints in this work, an interesting direction for future research is to investigate how input and state constraints can be satisfied by an OCO algorithm. Moreover, Assumption \ref{assump:TrackingCostFcn} could be relaxed, allowing economic cost functions, and more efficient algorithms than OGD could be applied. Finally, predictions on the future cost functions could be incorporated to improve the algorithm's performance.

\begin{small}
\bibliography{ifacconf}             % bib file to produce the bibliography

\begin{thebibliography}{4}
\providecommand{\natexlab}[1]{#1}
\providecommand{\url}[1]{\texttt{#1}}
\providecommand{\urlprefix}{URL }
\expandafter\ifx\csname urlstyle\endcsname\relax
  \providecommand{\doi}[1]{doi:\discretionary{}{}{}#1}\else
  \providecommand{\doi}{doi:\discretionary{}{}{}\begingroup
  \urlstyle{rm}\Url}\fi

\bibitem[{Able(1956)}]{Abl:56}
Able, B. (1956).
\newblock Nucleic acid content of microscope.
\newblock \emph{Nature}, 135, 7--9.

\bibitem[{Able et~al.(1954)Able, Tagg, and Rush}]{AbTaRu:54}
Able, B., Tagg, R., and Rush, M. (1954).
\newblock Enzyme-catalyzed cellular transanimations.
\newblock In A.~Round (ed.), \emph{Advances in Enzymology}, volume~2, 125--247.
  Academic Press, New York, 3rd edition.

\bibitem[{Keohane(1958)}]{Keo:58}
Keohane, R. (1958).
\newblock \emph{Power and Interdependence: World Politics in Transitions}.
\newblock Little, Brown \& Co., Boston.

\bibitem[{Powers(1985)}]{Pow:85}
Powers, T. (1985).
\newblock Is there a way out?
\newblock \emph{Harpers}, 35--47.

\end{thebibliography}


\begin{thebibliography}{24}
\providecommand{\natexlab}[1]{#1}
\providecommand{\url}[1]{\texttt{#1}}
\providecommand{\urlprefix}{URL }
\expandafter\ifx\csname urlstyle\endcsname\relax
  \providecommand{\doi}[1]{doi:\discretionary{}{}{}#1}\else
  \providecommand{\doi}{doi:\discretionary{}{}{}\begingroup
  \urlstyle{rm}\Url}\fi

\bibitem[{Abbasi-Yadkori et~al.(2014)Abbasi-Yadkori, Bartlett, and
  Kanade}]{Abbasi2014}
Abbasi-Yadkori, Y., Bartlett, P., and Kanade, V. (2014).
\newblock Tracking adversarial targets.
\newblock In \emph{Proceedings of the 31st International Conference on Machine
  Learning}, volume~32, 369--377.

\bibitem[{Agarwal et~al.(2019)Agarwal, Bullins, Hazan, Kakade, and
  Singh}]{Agarwal2019}
Agarwal, N., Bullins, B., Hazan, E., Kakade, S., and Singh, K. (2019).
\newblock Online control with adversarial disturbances.
\newblock In \emph{Proceedings of the 36th International Conference on Machine
  Learning}, volume~97, 111--119.

\bibitem[{{Akbari} et~al.(2019){Akbari}, {Gharesifard}, and
  {Linder}}]{Akbari2019}
{Akbari}, M., {Gharesifard}, B., and {Linder}, T. (2019).
\newblock {An Iterative Riccati Algorithm for Online Linear Quadratic Control}.
\newblock \emph{arXiv e-prints}.
\newblock ArXiv:1912.09451.

\bibitem[{Besbes et~al.(2015)Besbes, Gur, and Zeevi}]{Besbes2015}
Besbes, O., Gur, Y., and Zeevi, A. (2015).
\newblock Non-stationary stochastic optimization.
\newblock \emph{Operations Research}, 63(5), 1227--1244.

\bibitem[{Boyd and Vandenberghe(2004)}]{Boyd2004}
Boyd, S. and Vandenberghe, L. (2004).
\newblock \emph{Convex Optimization}.
\newblock Cambridge University Press, New York, NY, USA.

\bibitem[{Cao and Liu(2019)}]{Cao2019}
Cao, X. and Liu, K.J.R. (2019).
\newblock Online convex optimization with time-varying constraints and bandit
  feedback.
\newblock \emph{IEEE Transactions on automatic control}, 64(7), 2665--2680.

\bibitem[{Cohen et~al.(2018)Cohen, Hasidim, Koren, Lazic, Mansour, and
  Talwar}]{Cohen2018}
Cohen, A., Hasidim, A., Koren, T., Lazic, N., Mansour, Y., and Talwar, K.
  (2018).
\newblock Online linear quadratic control.
\newblock In \emph{Proceedings of the 35th International Conference on Machine
  Learning}, volume~80, 1029--1038.

\bibitem[{{Colombino} et~al.(2020){Colombino}, {Dall’Anese}, and
  {Bernstein}}]{Colombino2020}
{Colombino}, M., {Dall’Anese}, E., and {Bernstein}, A. (2020).
\newblock Online optimization as a feedback controller: Stability and tracking.
\newblock \emph{IEEE Transactions on Control of Network Systems}, 7(1),
  422--432.

\bibitem[{Diehl et~al.(2005)Diehl, Findeisen, Allg\"ower, Bock, and
  Schl\"oder}]{Diehl2005}
Diehl, M., Findeisen, R., Allg\"ower, F., Bock, H.G., and Schl\"oder, J.P.
  (2005).
\newblock Nominal stability of real-time iteration scheme for nonlinear model
  predictive control.
\newblock \emph{IEE Proceedings - Control Theory and Applications}, 152(3), 296
  -- 308.

\bibitem[{Faulwasser et~al.(2018)Faulwasser, Gr\"une, and
  M\"uller}]{Faulwasser2018}
Faulwasser, T., Gr\"une, L., and M\"uller, M.A. (2018).
\newblock Economic nonlinear model predictive control.
\newblock \emph{Foundations and Trends\textsuperscript{\textregistered} in
  Systems and Control}, 5(1), 1--98.

\bibitem[{Hall and Willett(2015)}]{Hall2015}
Hall, E.C. and Willett, R.M. (2015).
\newblock Online convex optimization in dynamic environments.
\newblock \emph{IEEE Journal of Selected Topics in Signal Processing}, 9(4),
  647--662.

\bibitem[{Hazan(2016)}]{Hazan2016}
Hazan, E. (2016).
\newblock Introduction to online convex optimization.
\newblock \emph{Foundations and Trends\textsuperscript{\textregistered} in
  Optimization}, 2(3-4), 157--325.

\bibitem[{Hazan et~al.(2007)Hazan, Agarwal, and Kale}]{Hazan2007}
Hazan, E., Agarwal, A., and Kale, S. (2007).
\newblock Logarithmic regret algorithms for online convex optimization.
\newblock \emph{Machine Learning}, 69(2), 169--192.

\bibitem[{Li et~al.(2019)Li, Chen, and Li}]{Li2019}
Li, Y., Chen, X., and Li, N. (2019).
\newblock Online optimal control with linear dynamics and predictions:
  Algorithms and regret analysis.
\newblock In \emph{Advances in Neural Information Processing Systems 32},
  14887--14899. Curran Associates, Inc.

\bibitem[{Li et~al.(2018)Li, Qu, and Li}]{Li2018}
Li, Y., Qu, G., and Li, N. (2018).
\newblock Using predictions in online optimization with switching costs: A fast
  algorithm and a fundamental limit.
\newblock \emph{2018 Annual American Control Conference (ACC)}, 3008--3013.

\bibitem[{Lin et~al.(2013)Lin, Wierman, Andrew, and Thereska}]{Lin2013}
Lin, M., Wierman, A., Andrew, L.L.H., and Thereska, E. (2013).
\newblock Dynamic right-sizing for power-proportional data centers.
\newblock \emph{IEEE/ACM Transactions on Networking}, 21(5), 1378--1391.

\bibitem[{Mokhtari et~al.(2016)Mokhtari, Shahrampour, Jadbabaie, and
  Ribeiro}]{Mokhtari2016}
Mokhtari, A., Shahrampour, S., Jadbabaie, A., and Ribeiro, A. (2016).
\newblock Online optimization in dynamic environments: Improved regret rates
  for strongly convex problems.
\newblock \emph{2016 IEEE 55th Conference on Decision and Control (CDC)},
  7195--7201.

\bibitem[{Nesterov(2018)}]{Nesterov2018}
Nesterov, Y. (2018).
\newblock \emph{Lectures on Convex Optimization}, volume 137 of \emph{Springer
  Optimization and Its Applications}.
\newblock Springer International Publishing, 2 edition.

\bibitem[{Paternain and Ribeiro(2016)}]{Paternain2016}
Paternain, S. and Ribeiro, A. (2016).
\newblock Online learning of feasible strategies in unknown environments.
\newblock \emph{IEEE Transactions on Automatic Control}, 62(6), 2807--2822.

\bibitem[{Rawlings and Mayne(2009)}]{Rawlings2009}
Rawlings, J.B. and Mayne, D.Q. (2009).
\newblock \emph{Model Predictive Control: Theory and Design}.
\newblock Nob Hill Pub.

\bibitem[{Scokaert et~al.(1999)Scokaert, Mayne, and Rawlings}]{Scokaert1999}
Scokaert, P.O.M., Mayne, D.Q., and Rawlings, J.B. (1999).
\newblock Suboptimal model predictive control (feasibility implies stability).
\newblock \emph{IEEE Transactions on Automatic Control}, 44(3), 648--654.

\bibitem[{Shalev-Shwartz(2012)}]{Shwartz2012}
Shalev-Shwartz, S. (2012).
\newblock Online learning and online convex optimization.
\newblock \emph{Foundations and Trends\textsuperscript{\textregistered} in
  Machine Learning}, 4(2), 107--194.

\bibitem[{Tanaka(2006)}]{Tanaka2006}
Tanaka, M. (2006).
\newblock Real-time pricing with ramping costs: A new approach to managing a
  steep change in electricity demand.
\newblock \emph{Energy Policy}, 34(18), 3634--3643.

\bibitem[{Zinkevich(2003)}]{Zinkevich2003}
Zinkevich, M. (2003).
\newblock Online convex programming and generalized infinitesimal gradient
  ascent.
\newblock \emph{Proceedings of the Twentieth International Conference on
  Machine Learning (ICML)}, 928 -- 936.

\end{thebibliography}
\end{small}     

\appendix
\section{Proof of Theorem \ref{thm}} \label{appendix}

Before the formal proof of Theorem \ref{thm}, we derive three auxiliary results. First, we study the closed-loop dynamics resulting from application of Algorithm~1 to system~\eqref{eq:sysdyn}. Let $t \in \mathbb{N}_{[1,T-\mu+1]}$, then we have by repeatedly applying the system dynamics \eqref{eq:sysdyn}
\begin{align*}
x_{t+\mu-1} &= A^\mu x_{t-1} + \sum_{i=0}^{\mu-1} A^iBu_{t+\mu-1-i} \\
&\refeq{\eqref{eq:GenB_Output}} A^\mu x_{t-1} + S_c V_t^\mu + Be g_{t+\mu-1} \\ &\hspace{15pt}+ (BeW+ABe) g_{t+\mu-2} + \dots \\
&\hspace{15pt}+ (BeW^{\mu-1} + ABeW^{\mu-2} + \dots + A^{\mu-1}Be)g_t \\
&\hspace{15pt}+ \dots + A^{\mu-1}BeW^{\mu-1}g_{t-\mu+1},
\end{align*}
where $V^\mu_t = \begin{pmatrix} v_{t+\mu-1}^T & \dots & v_t^T \end{pmatrix}^T \in \mathbb{R}^{\mu m}$. 

Inserting the relations $\sum_{i=0}^k A^iBeW^{k-i} = S_c(W^T)^{\mu-1-k}$ and $\sum_{i=0}^k A^{\mu-1-i}BeW^{\mu-1-k+i} = S_cW^{\mu-1-k}$, where \linebreak $k\in\mathbb{N}_{[0,\mu-1]}$, yields
\begin{align}
\begin{split}
x_{t+\mu-1} = &A^\mu x_{t-1} + S_c V_t^\mu + S_c g_t \\ &+ S_c \sum_{i=1}^{\mu-1}(W^T)^ig_{t+i} + S_c \sum_{i=1}^{\mu-1}W^i g_{t-i}. \label{eq:cldynamics}
\end{split}
\end{align}

Next, the predictions in \eqref{eq:GenB_StatePrediction} can be calculated recursively. Let $E_{01} = \begin{pmatrix} \bm 0_{m\times(\mu-1)m} \\ \bm I_{(\mu-1)m} \end{pmatrix}$, then we have for $t\in\mathbb{N}_{[\mu,T+\mu-2]}$
\begin{align*}
\hat x_{t+1} &\refeq{\eqref{eq:GenB_StatePrediction}} A^\mu x_{t-\mu+1}+S_c\hat V_{t-\mu+2} + S_c \sum_{i=1}^{\mu-1} W^i g_{t-\mu+2-i} \\
&\refeq{\eqref{eq:sysdyn},\eqref{eq:GenB_Output}} \hspace{2pt} A^\mu \left(Ax_{t-\mu}+Bv_{t-\mu+1}+B\sum_{i=0}^{\mu-1}eW^ig_{t-\mu+1-i} \right) \\
&\hspace{10pt}+ \sum_{i=0}^{\mu-1}A^iBv_{t-\mu+2} + S_cW\sum_{i=0}^{\mu-1} W^ig_{t-\mu+1-i} \\
&= A \left( A^\mu x_{t-\mu} + S_c \hat V_{t-\mu+1} + S_c \sum_{i=1}^{\mu-1} W^ig_{t-\mu+1-i} \right) \\
&\hspace{10pt} - AS_c \hat V_{t-\mu+1} + A^\mu Bv_{t-\mu+1} + \sum_{i=0}^{\mu-1} A^iBv_{t-\mu+2} \\
&\hspace{10pt}+ AS_c g_{t-\mu+1},
\end{align*}
where we use $W^\mu = 0$ in the second and the relation $S_cW+A^\mu Be = AS_c$ in the third line. Inserting \eqref{eq:GenB_StatePrediction} yields
\begin{align}
\begin{split}
\hat x_{t+1} = &A \hat x_t + Bv_{t-\mu+2} + AS_cg_{t-\mu+1} \\
&+ S_cE_{01}E_{01}^T \left(\hat V_{t-\mu+2} - \hat V_{t-\mu+1} \right).
\end{split} \label{eq:recpred}
\end{align}

\vspace{-1ex}
Third, we have the following result on the rate of convergence of gradient descent \citep[Chapter 2.3.3, Theorem 2.3.4]{Nesterov2018}. For an $\alpha$-convex and $l$-smooth function $f: \mathbb{R}^n \rightarrow \mathbb{R}$ to be minimized, one gradient descent step $x_1 = x_0 - \gamma \nabla f(x_0)$ yields
\begin{align}
\norm{x_1 - \theta} \leq \kappa \norm{x_0 - \theta}, \label{eq:Nesterov}
\end{align}
where $\theta = \arg \min_x f(x)$ and $\kappa = 1-\alpha\gamma$. Accordingly, we define $\kappa_x = 1-\alpha_x\gamma_x$ and $\kappa_v = 1-\alpha_u \gamma_v$. Before we prove Theorem \ref{thm}, we introduce the following supporting lemma.
\begin{lem} \label{lemma}
	Let Assumptions \ref{assump:CostFunctionProp}, \ref{assump:TrackingCostFcn}, and \ref{assump:SysCond} be satisfied. Let $\alpha_u > \frac{2-\sqrt 2}{2+\sqrt2} l_u$. Given step sizes $\frac{2-\sqrt 2}{2\alpha_u} < \gamma_v \leq \frac{2}{l_u+\alpha_u}$ and $\frac{2\norm{A}-1}{2\norm{A}\alpha_x} < \gamma_x \leq \frac{2}{l_x + \alpha_x}$, the predicted states $\hat x_t$ satisfy
	\begin{align*}
	\sum_{t=1}^\tau \norm{\hat x_{t+\mu-1}-\theta_{t-1}}^2 \leq &C_\theta/l_x \sum_{t=1}^{\tau-1} \norm{\theta_t - \theta_{t-1}}^2 \\ &+ C_\eta/l_x \sum_{t=1}^{\tau-1} \norm{\eta_t - \eta_{t-1}}^2,
	\end{align*}
	where $\tau \in \mathbb{N}_{[1,T]}$, $C_\eta = \frac{8l_x(\norm{B}^2\kappa_v^2 + \norm{S_cE_{01}}^2 \gamma_v^2l_u^2(\mu-1))}{(1-2\kappa_v^2) (1-4\norm{A}^2\kappa_x^2)}$ and $C_\theta = \frac{4\norm{A}^2l_x}{1-4\norm{A}^2\kappa_x^2}$.
\end{lem}

\emph{Proof.} Note that the bounds on the step sizes imply $1-2\kappa_v^2 > 0$ and $1-4\norm{A}^2\kappa_x^2 >0$. By Jensen's inequality we have
\begin{align}
\norm{\sum_{i=1}^p a_i}^2 \leq \sum_{i=1}^p \lambda_i \norm{a_i}^2 \label{eq:triangle},
\end{align}
where $\sum_{i=1}^p \frac{1}{\lambda_i} = 1$ and $\lambda_i > 0$ for all $i\in\mathbb{N}_{[1,p]}$. In particular, we can choose $\lambda_i = p$ for all $i \in \mathbb{N}_{[1,p]}$.

Next, we derive three auxiliary results on the relation of $v_t$ and $\eta_t$. Remember that $\eta_0 = v_0 = v_1$. Hence, we have
\begin{align*}
&\sum_{t=1}^{\tau-1} \norm{v_{t+1}-\eta_t}^2 \refeq{\eqref{eq:GenB_InputOGD}} \sum_{t=1}^{\tau-1} \norm{v_{t} - \gamma_v \nabla f_t^u(v_t) - \eta_t}^2 \\
&\refleq{\eqref{eq:Nesterov}} \hspace{5pt} \kappa_v^2 \sum_{t=1}^{\tau-1} \norm{v_t - \eta_t}^2 \\
&\refleq{\eqref{eq:triangle}} \hspace{5pt} 2\kappa_v^2 \sum_{t=1}^{\tau-1} \norm{v_t - \eta_{t-1}}^2 + 2\kappa_v^2 \sum_{t=1}^{\tau-1} \norm{\eta_t - \eta_{t-1}}^2 \\
&\refleq{v_1 = \eta_0} \hspace{5pt} 2\kappa_v^2 \sum_{t=1}^{\tau-1} \norm{v_{t+1} - \eta_t}^2 + 2\kappa_v^2 \sum_{t=1}^{\tau-1} \norm{\eta_t - \eta_{t-1}}^2.
\end{align*}
Since $1-2\kappa_v^2 > 0$, rearranging yields
\begin{align}
\sum_{t=1}^{\tau-1} \norm{v_{t+1} - \eta_t}^2 &\leq \frac{2\kappa_v^2}{1-2\kappa_v^2} \sum_{t=1}^{\tau-1} \norm{\eta_t - \eta_{t-1}}^2. \label{eq:v+1-eta}
\end{align}
Moreover, we have
\begin{align*}
\sum_{t=1}^\tau \norm{v_t - \eta_t}^2 \hspace{5pt} &\refleq{\eqref{eq:triangle}} \hspace{7pt} 2\sum_{t=1}^\tau \norm{v_t - \eta_{t-1}}^2 + 2\sum_{t=1}^\tau \norm{\eta_t - \eta_{t-1}}^2 \\
&\refeq{v_1 = \eta_0} \hspace{7pt} 2\sum_{t=1}^{\tau-1} \norm{v_{t+1}-\eta_t}^2 + 2\sum_{t=1}^\tau\norm{\eta_t - \eta_{t-1}}^2 \\
&\refleq{\eqref{eq:v+1-eta}} \hspace{7pt} \frac{2}{1-2\kappa_v^2} \sum_{t=1}^\tau \norm{\eta_t - \eta_{t-1}}^2. \numberthis \label{eq:vt-etat}
\end{align*}
Finally, due to optimality of $\eta_t$ and, hence, $\nabla f_t^u(\eta_t) = 0$,
\begin{align*}
&\sum_{t=1}^{k-1} \norm{v_{t+1}-v_t}^2 \refleq{\eqref{eq:GenB_InputOGD}} \gamma_v^2 \sum_{t=1}^{k-1} \norm{\nabla f_t^u(v_t) - \nabla f_t^u(\eta_t)}^2 \\
\leq &\gamma_v^2 l_u^2 \sum_{t=1}^{k-1} \norm{v_t - \eta_t}^2 \hspace{1ex} \refleq{\eqref{eq:vt-etat}} \frac{2 \gamma_v^2 l_u^2}{1-2\kappa_v^2} \sum_{t=1}^{k-1} \norm{\eta_t - \eta_{t-1}}^2, \numberthis \label{eq:v+1-v}
\end{align*}
where we use the fact that $l_u$-smoothness of $f_t^u(u)$ implies a Lipschitz condition on its gradient in the second line.

Last, we show the bound on the predicted states. Since $\theta_0 = \hat x_\mu$, we have
\begingroup
\allowdisplaybreaks
\begin{align*}
\sum_{t=1}^\tau &\norm{\hat x_{t+\mu-1} - \theta_{t-1}}^2 = \sum_{t=1}^{\tau-1} \norm{\hat x_{t+\mu}-\theta_t} \\ 
\begin{split}
&\refleq{\eqref{eq:recpred},\eqref{eq:triangle}} \hspace{15pt} \sum_{t=1}^{\tau-1} \bigg( 2\norm{A(\hat x_{t+\mu-1} + S_c g_t - \theta_t)}^2 \bigg. \\ &\hspace{5pt} \left. + 4\norm{B(v_{t+1} - \eta_t)}^2 + 4\norm{S_c E_{01} E_{01}^T \left(\hat V_{t+1} - \hat V_t \right)}^2 \right) \end{split} \\
&\refleq{\eqref{eq:triangle},\eqref{eq:v+1-eta}} \hspace{12pt} 4\norm{A}^2 \sum_{t=1}^{\tau-1} \norm{\hat x_{t+\mu-1} - \gamma_x \nabla f_{t-1}^x(\hat x_{t+\mu-1}) - \theta_{t-1}}^2 \\ &\hspace{2pt}+ \hspace{-2pt} 4\norm{A}^2 \sum_{t=1}^{\tau-1} \norm{\theta_t - \theta_{t-1}}^2 \hspace{-3pt} + \hspace{-2pt} \frac{8 \norm{B}^2 \kappa_v^2}{1-2\kappa_v^2} \sum_{t=1}^{\tau-1} \norm{\eta_t-\eta_{t-1}}^2 \\ &\hspace{2pt} + 4 \norm{S_cE_{01}}^2 (\mu-1) \sum_{t=1}^{\tau-1} \norm{v_{t+1} - v_t}^2 \\ 
\begin{split}
&\refleq{\eqref{eq:Nesterov},\eqref{eq:v+1-v}} \hspace{12pt} 4 \norm{A}^2 \kappa_x^2 \sum_{t=1}^{\tau-1} \norm{\hat x_{t+\mu-1} - \theta_{t-1}}^2 \\ &\hspace{20pt} + 4\norm{A}^2 \sum_{t=1}^{\tau-1} \norm{\theta_t - \theta_{t-1}}^2 \\ &\hspace{20pt}+ \frac{C_\eta (1-4\norm{A}^2\kappa_x^2)}{l_x} \sum_{t=1}^{\tau-1} \norm{\eta_t - \eta_{t-1}}^2, \end{split}
\end{align*}

\vspace{-2ex}

Since $1-4\norm{A}^2\kappa_x^2 >0$, rearranging yields the desired bound. \hfill\hfill\qed
\endgroup

\emph{Proof of Theorem \ref{thm}.} The proof consists of three parts. First, we derive a bound on the cost of the control inputs. In the second part, we derive a bound on $\sum_{t=1}^k \norm{\hat x_{t+\mu-1} - x_{t+\mu-1}}^2$, which will be useful to bound the cost on the states. The last part is to combine these results to find a bound on the regret of Algorithm~1.

First, we have for $k \in \mathbb{N}_{[1,T]}$
\begin{align*}
&\sum_{t=1}^k \norm{\sum_{i=0}^{\mu-1} eW^ig_{t-i}}^2 \refleq{\eqref{eq:triangle}} \mu \sum_{t=1}^k \sum_{i=0}^{\mu-1} \norm{eW^i g_{t-i}}^2 \\
\leq &\mu \sum_{t=1}^k \sum_{i=0}^{\mu-1} \norm{eW^ig_t}^2 \\
\refleq{\eqref{eq:GenB_StateOGD}} &\mu\gamma_x^2C_1 \sum_{t=1}^k \norm{\nabla f_{t-1}^x(\hat x_{t+\mu-1}) - \nabla f_{t-1}^x(\theta_{t-1})}^2 \\
\leq &\mu\gamma_x^2l_x^2 C_1 \sum_{t=1}^k \norm{\hat x_{t+\mu-1} - \theta_{t-1}}^2,
\end{align*}
where we use the fact that $g_t=0$ for $t\leq1$ in the second line, $C_1 = \sum_{i=0}^{\mu-1} \norm{eW^iS_c^T(S_cS_c^T)^{-1}}^2$ and \mbox{$\nabla f_{t-1}^x(\theta_{t-1})=0$} in the third line, and Lipschitz continuity of the gradients in the last line. By Lemma \ref{lemma}, we obtain
\begin{align}
\begin{split}
\sum_{t=1}^k \norm{\sum_{i=0}^{\mu-1} eW^ig_{t-i}}^2 &\leq \mu \gamma_x^2 l_xC_1 C_\theta \sum_{t=1}^{k-1} \norm{\theta_t - \theta_{t-1}}^2 \\ &\hspace{20pt} + \mu\gamma_x^2l_xC_1C_\eta \sum_{t=1}^{k-1} \norm{\eta_t-\eta_{t-1}}^2.
\end{split} \label{eq:sumewigt-i}
\end{align}
\vspace{-1ex}
Hence, we have
\begin{align*}
&\sum_{t=1}^T \hspace{-2pt} \norm{u_t - \eta_t}^2 \hspace{10pt} \refleq{\eqref{eq:GenB_Output},\eqref{eq:triangle}} \hspace{8pt} 2\sum_{t=1}^T \hspace{-2pt} \norm{v_t - \eta_t}^2 \hspace{-3pt} + \hspace{-2pt} 2 \hspace{-2pt} \sum_{t=1}^T \norm{\sum_{i=0}^{\mu-1} \hspace{-2pt} eW^i \hspace{-2pt} g_{t-i}}^2 \\
\begin{split}
& \hspace{15pt} \refleq{\eqref{eq:vt-etat}, \eqref{eq:sumewigt-i}} \hspace{15pt} \frac{4}{1-2\kappa_v^2}\sum_{t=1}^T \norm{\eta_t - \eta_{t-1}}^2 \\ &+ \hspace{-2pt} 2\mu\gamma_x^2l_xC_1 \hspace{-2pt} \left( \hspace{-2pt} C_\theta \sum_{t=1}^{T-1} \norm{\theta_t - \theta_{t-1}}^2 + C_\eta \sum_{t=1}^{T-1} \norm{\eta_t-\eta_{t-1}}^2 \right)\hspace{-3pt}.
\end{split} \numberthis \label{eq:ut-etat2}
\end{align*}

Having established a bound on the regret of the control inputs chosen by Algorithm~1, we derive a bound on $\sum_{t=1}^k \norm{\hat x_{t+\mu-1} - x_{t+\mu-1}}^2$ in the second part of the proof. First, let $k\in\mathbb{N}_{[1,T-\mu+1]}$. The last component of $\hat V_t$ and $V_t^\mu$ is the same, therefore, we insert the matrix \mbox{$E_{10} = \begin{pmatrix} I_{(\mu-1)m} \\ \bm{0}_{m \times (\mu-1)m} \end{pmatrix} \in \mathbb{R}^{m \mu \times (\mu-1)m}$} which yields
\vspace{-1ex}
\begin{align*} % Das hier ist lang - unterbrechen und Erklärung einfügen?
&\sum_{t=1}^k \norm{S_c\left(\hat V_t - V_t^\mu \right)}^2 = \sum_{t=1}^k \norm{S_cE_{10}E_{10}^T \left(\hat V_t - V_t^\mu \right)}^2 \\
&\leq \norm{S_cE_{10}}^2 \sum_{t=1}^k \sum_{i=1}^{\mu-1} \norm{v_{t+i} - v_t}^2 \\
&= \norm{S_cE_{10}}^2 \sum_{t=1}^k \sum_{i=1}^{\mu-1} \norm{\sum_{j=1}^{i} v_{t+j} - v_{t+j-1}}^2,
\end{align*}
where we make use of a telescoping series in the last line. Next, by inserting \eqref{eq:triangle}, upper bounding $i$, and finally positivity of the norm we have
\begin{align*}
&\sum_{t=1}^k \norm{S_c\left(\hat V_t - V_t^\mu \right)}^2 \\
& \hspace{20pt} \refleq{\eqref{eq:triangle}} \hspace{3pt} \norm{S_cE_{10}}^2 \sum_{t=1}^k \sum_{i=1}^{\mu-1} i \sum_{j=1}^{i} \norm{v_{t+j} - v_{t+j-1}}^2 \\ % why this term and not ||S_cE_{10}E_{10}^T||
& \hspace{20pt} \leq \norm{S_cE_{10}}^2 (\mu-1)^2\sum_{t=1}^k\sum_{j=1}^{\mu-1} \norm{v_{t+j} - v_{t+j-1}}^2 \\
& \hspace{20pt} \leq \norm{S_cE_{10}}^2 (\mu-1)^3\sum_{t=1}^{k+\mu-2} \norm{v_{t+1} - v_{t}}^2 \\
& \hspace{20pt} \refleq{\eqref{eq:v+1-v}} \hspace{3pt} C_2 \hspace{-2pt} \sum_{t=1}^{k+\mu-2} \norm{\eta_t - \eta_{t-1}}^2, \numberthis \label{eq:SchatV-Vtmu}
\end{align*}
where $C_2 = \frac{2\norm{S_cE_{10}}^2\gamma_v^2l_u^2(\mu-1)^3}{1-2\kappa_v^2}$. Next, we apply \eqref{eq:triangle} to obtain
\begin{align*}
&\sum_{t=1}^k \norm{S_c \sum_{i=0}^{\mu-1} (W^T)^i g_{t+i}}^2 \refleq{\eqref{eq:triangle}} \sum_{t=1}^k \sum_{i=0}^{\mu-1} \mu \norm{S_c(W^T)^ig_{t+i}}^2 \\
\leq &\mu \sum_{t=1}^{k+\mu-1} \sum_{i=0}^{\mu-1} \norm{S_c(W^T)^ig_{t}}^2 \\
\refleq{\eqref{eq:GenB_StateOGD}} &\mu\gamma_x^2 C_3 \sum_{t=1}^{k+\mu-1} \norm{\nabla f_{t-1}^x(\hat x_{t+\mu-1}) - \nabla f_{t-1}^x(\theta_{t-1})}^2,
\end{align*}
where $C_3 = \sum_{i=0}^{\mu-1} \norm{S_c(W^T)^iS_c^T(S_cS_c^T)^{-1}}^2$ and we use the fact that $\nabla f_t^x(\theta_t) = 0$ due to optimality of $\theta_t$ in the last line. Lipschitz continuity of the gradient yields
\begin{align*}
\sum_{t=1}^k \norm{S_c \sum_{i=0}^{\mu-1} (W^T)^i g_{t+i}}^2\hspace{-5pt} &\leq \mu \gamma_x^2 l_x^2 C_3 \hspace{-5pt} \sum_{t=1}^{k+\mu-1} \norm{\hat x_{t+\mu-1} - \theta_{t-1}}^2 \hspace{-5pt}.
\end{align*}
By Lemma \ref{lemma}, we obtain
\begin{align}
\begin{split}
\sum_{t=1}^k \norm{S_c \sum_{i=0}^{\mu-1} (W^T)^i g_{t+i}}^2 &\leq \mu\gamma_x^2l_xC_3C_\theta \sum_{t=1}^{k+\mu-2} \norm{\theta_t - \theta_{t-1}}^2 \\ &+ \mu\gamma_x^2l_xC_3C_\eta \sum_{t=1}^{k+\mu-2} \norm{\eta_t - \eta_{t-1}}^2. 
\end{split} \label{eq:ScsumWTigt+i}
\end{align}
Combining the last two results yields the desired bound
\begin{align*}
\sum_{t=1}^k &\norm{\hat x_{t+\mu-1} - x_{t+\mu-1}}^2 \\ &\refeq{\eqref{eq:GenB_StatePrediction},\eqref{eq:cldynamics}} \hspace{10pt} \sum_{t=1}^{k} \norm{S_c\left(\hat V_t - V_t^\mu \right) - S_c \sum_{i=0}^{\mu-1} (W^T)^i g_{t+i}}^2 \\
&\refleq{\eqref{eq:triangle}} \hspace{2pt} 2 \sum_{t=1}^k \norm{S_c \left( \hat V_t - V_t^\mu \right) }^2 \hspace{-5pt} + \hspace{-2pt} 2 \hspace{-2pt} \sum_{t=1}^k \norm{S_c \sum_{i=0}^{\mu-1} (W^T)^i g_{t+i}}^2 \\
\begin{split}
&\refleq{\eqref{eq:SchatV-Vtmu},\eqref{eq:ScsumWTigt+i}} \hspace{18pt} 2\mu\gamma_x^2l_xC_3C_\theta \sum_{t=1}^{k+\mu-2} \norm{\theta_t - \theta_{t-1}}^2 \\ &\hspace{5pt}+ 2\left( C_2 + \mu\gamma_x^2l_xC_3C_\eta \right) \sum_{t=1}^{k+\mu-2} \norm{\eta_t - \eta_{t-1}}^2.
\end{split} \numberthis \label{eq:hatxt+mu-1-xt+mu-1}
\end{align*}
Before we combine all results to obtain a bound on the regret, we apply \eqref{eq:triangle} to obtain
\begin{align*}
\sum_{t=1}^k \norm{\theta_{t+p}-\theta_{t-1}}^2 &= \hspace{5pt} \sum_{t=1}^k \norm{\sum_{i=0}^p \theta_{t+i} - \theta_{t+i-1}}^2 \\
&\refleq{\eqref{eq:triangle}} \hspace{5pt} \sum_{t=1}^k (p+1) \sum_{i=0}^p \norm{\theta_{t+i}-\theta_{t+i-1}}^2 \\
&\leq (p+1)^2 \sum_{t=1}^{k+p} \hspace{5pt} \norm{\theta_t - \theta_{t-1}}^2. \numberthis \label{eq:JensenTheta2}
\end{align*}
Finally, we are ready to compute the regret of Algorithm~1. By optimality of $(\theta_t,\eta_t)$, we have
\begin{align*}
\mathcal{R} &= \sum_{t=1}^T f_t^x(x_t) - f_t^x(x^*_t) + f_t^u(u_t) - f_t^u(u^*_t) \\
&\leq \sum_{t=1}^T f_t^x(x_t) - f_t^x(\theta_t) + f_t^u(u_t) - f_t^u(\eta_t).
\end{align*}
Next, we apply $l$-smoothness of the cost functions $f_t^x(x)$ and $f_t^u(u)$. Due to $\nabla f_t^x(\theta_t) = 0$, $\nabla f_t^u(\eta_t) = 0$, and after splitting up the sums we get
\begin{align*}
\mathcal{R} &\leq \frac{l_x}{2} \hspace{-2pt} \sum_{t=1}^{\mu-1} \norm{x_t - \theta_t}^2 \hspace{-3pt} + \hspace{-3pt}  \frac{l_x}{2} \hspace{-2pt} \sum_{t=\mu}^T \norm{x_t - \theta_t}^2 \hspace{-3pt} + \hspace{-3pt}  \frac{l_u}{2} \sum_{t=1}^T \norm{u_t - \eta_t}^2 \\
&= C_\mu \hspace{-2pt} + \hspace{-2pt} \frac{l_x}{2} \sum_{t=1}^{T-\mu+1} \hspace{-3pt} \norm{x_{t+\mu-1} - \theta_{t+\mu-1}}^2 \hspace{-2pt} + \hspace{-2pt} \frac{l_u}{2} \sum_{t=1}^T \norm{u_t - \eta_t}^2 \hspace{-3pt} .
\end{align*}
We first apply \eqref{eq:triangle} to the first sum and \eqref{eq:ut-etat2} to bound the cost on the control inputs. Afterwards, inserting \eqref{eq:triangle} and \eqref{eq:hatxt+mu-1-xt+mu-1} yields
\begin{align*}
\mathcal{R} \hspace{15pt} &\refleq{\eqref{eq:triangle},\eqref{eq:ut-etat2}} \hspace{15pt} C_\mu + l_x \sum_{t=1}^{T-\mu+1} \norm{\hat x_{t+\mu-1} - \theta_{t+\mu-1}}^2 \\ &\hspace{20pt}+ l_x \sum_{t=1}^{T-\mu+1} \norm{\hat x_{t+\mu-1} - x_{t+\mu-1}}^2 + \frac{2l_u}{1-2\kappa_v^2} H_T \\ &\hspace{20pt}+ \mu \gamma_x^2 l_x l_u C_1 \left( C_\theta \Theta_{T-1} + C_\eta H_{T-1} \right) \\
&\refleq{\eqref{eq:triangle},\eqref{eq:hatxt+mu-1-xt+mu-1}} \hspace{15pt} C_\mu + 2l_x \sum_{t=1}^{T-\mu+1} \norm{\hat x_{t+\mu-1} - \theta_{t-1}}^2 \\ &\hspace{20pt}+ 2l_x \sum_{t=1}^{T-\mu+1} \norm{\theta_{t+\mu-1}-\theta_{t-1}}^2 + C_4C_\theta \Theta_{T-1} \\ &\hspace{20pt}+ (2l_xC_2 + C_4C_\eta) H_{T-1} + \frac{2l_u}{1-2\kappa_v^2} H_T,
\end{align*}
where $C_4 = \mu\gamma_x^2l_x(2l_xC_3 + l_uC_1)$. Finally, we can apply Lemma \ref{lemma} and \eqref{eq:JensenTheta2} to obtain
\begin{align*}
\mathcal{R} \leq &C_\mu +2C_\theta \Theta_{T-\mu} + C_4C_\theta\Theta_{T-1} + 2l_x\mu^2 \Theta_T \\
&+ 2 C_\eta H_{T-\mu} + (2l_xC_2+C_4C_\eta) H_{T-1} + \frac{2l_u}{1-2\kappa_v^2}H_T.
\end{align*}
Positivity of the norm yields
\begin{align*}
	\mathcal{R} &\leq C_\mu + \Lambda_\theta \Theta_T + \Lambda_\eta H_T,
\end{align*}
where $\Lambda_\theta = 2C_\theta + C_4 C_\theta + 2l_x\mu^2$ and $\Lambda_\eta = 2C_\eta + 2l_xC_2 + C_4 C_\eta + \frac{2l_u}{1-2\kappa_v^2}$. \hfill \hfill \qed

\end{document}